\documentclass[preprint,12pt]{elsarticle} 
\usepackage[latin1]{inputenc}
\usepackage[english]{babel}
\usepackage{enumerate}           
\usepackage{amsmath}  
\usepackage{amssymb}   
\usepackage{graphicx}       
\usepackage{epsfig}         
\usepackage{color}      
\usepackage{caption}
\usepackage{float}
\usepackage{tabulary} 
\usepackage{csquotes}

\usepackage{enumerate}   

\journal{the scientific journal}   
\begin{document}

\makeatletter
\newcommand{\norm}[1]{\left\lVert#1\right\rVert}
\makeatother
  
%%%%%%%%%%%%%%%%%

\begin{abstract}
In this paper, we establish the sufficient conditions guaranteeing global uniform exponential stability, or at least global asymptotic stability, of all solutions for nonlinear dynamical systems, also known as global incremental stability (GIS) of the systems. We provide here an alternative approach for assessment of GIS in terms of logarithmic norm under which the stability becomes a topological notion and also generalize both horizontally and vertically the well-known Demidovich criterion for GIS of dynamical systems. Convergence of all solutions to the origin $x=0,$ which is not assumed to be an equilibrium state of system, is also analyzed. Theory is illustrated by a simulation experiment.
\end{abstract}  

\title{A novel criterion for global incremental stability of dynamical systems}
\author{Robert~Vrabel}
\ead{robert.vrabel@stuba.sk}
\address{Slovak University of Technology in Bratislava, Institute of Applied Informatics, Automation and Mechatronics,  Bottova 25, 917 01 Trnava, SLOVAKIA}

\newtheorem{thm}{Theorem}
\newtheorem{lem}{Lemma}
\newtheorem{defi}{Definition}
\newdefinition{rmk}{Remark}
\newdefinition{ex}{Example}
\newproof{pf}{Proof}
\newproof{proofoflemma2}{Proof of Lemma~\ref{integral_eq}}

\pagestyle{headings}

\maketitle
\section[Introduction]{Introduction}
Undoubtedly, stability analysis is one of the most important topics in dynamical systems theory. Traditionally, the stability of {\it a particular solution} of dynamical systems is analyzed, most often the origin, for example if we study the global error dynamics in the state trajectory tracking problem \cite{Liu_Huang}, \cite{Mazenc}. Stability analysis in the global sense could also have its benefits in the study of the convergence to zero of {\it all solution} in the cases, when the origin is not a solution of perturbed system, for example, in the situation when the origin is a stable equilibrium of nominal (unperturbed) system and we are interested in the effect of an external disturbance on the behavior of the systems as a part of the robustness analysis. Some new results in this field are direct consequences of the second part of Theorem~\ref{theorem_main} and demonstrated in Example~\ref{example1}. 

\section{ Notations and preliminaries} 
Our purpose here is to prove a new result regarding the global asymptotic and global uniform exponential stability of {\it all solutions} of perturbed nonlinear system 
\begin{equation}\label{original_system}
\dot x=f(x,t)+\delta(t),\quad x\in\mathbb{R}^n,  \quad t\geq t_0, 
\end{equation}
given that $x=0$ may not be a solution for the nominal system $\dot x=f(x,t)$ and that nominal vector field $f$ and  perturbation $\delta$ satisfy certain conditions described in the terms of logarithmic norm.  In other words, we focus here on the systems whose trajectories converge to one another and, in general, without being attracted toward some equilibrium position. The underlying idea is obvious: {\it If we have proved the asymptotic stability of all solutions at once, we do not have to deal with the stability properties of a particular solution}, especially if finding it itself is a difficult task \cite{Pavlov}. 

Two similar, but not entirely equivalent \cite{Ruffer} stability notions were settled - one is the long established notion of convergent systems \cite{Demidovich}, \cite{Pavlov}, the other is the younger notion of incremental stability \cite{Angeli}.  
\begin{defi}[cf.~\cite{Ruffer}] \label{GIS} 
System (\ref{original_system}) is incrementally asymptotically stable in a positively invariant set $X\subset\mathbb{R}^n$ if there exists a function $\beta\in{\cal K}{\cal L}$ \cite{Khalil} such that for any two solution $ x(t)$ and $x^*(t)$ with $x(t_0), x^*(t_0)\in X$ and $t\geq t_0,$
\[
\vert x(t)-x^*(t)\vert\leq\beta\big(\vert x(t_0)-x^*(t_0)\vert, t-t_0 \big).
\]
In the case $X = \mathbb{R}^n$ we say that system (\ref{original_system}) is globally incrementally stable (GIS).
\end{defi}
The aim of this paper is to provide an alternative approach for assessment of GIS based on the logarithmic norm and the variation of constant formula applied to auxiliary linear time-varying systems. We obtain more general results  as those achieved by using (quadratic) Lyapunov-like function which until now has been practically the only applicable method, see, e. g. \cite{Pavlov}, \cite{Ruffer} and the references therein. Moreover, the proposed approach turns out to be a bit simpler than through finding some implicit motion integral as in Lyapunov theory.

As a completely new result seems to be the establishing of conditions for the convergence of all solutions of the system (\ref{original_system}) to $0$ as $t\to\infty$ if even $x=0$ is not the equilibrium position of the nominal system $\dot x=f(x,t)$, in Theorem~\ref{theorem_main}; the context and novelty are explained in Remarks~\ref{Demidovich1} and~\ref{Demidovich2}. From another point of view, if $f(0,t)=0$ for all $t\geq t_0,$ Theorem~\ref{theorem_main} gives a sufficient conditions for robustness of global asymptotic stability to external perturbation $\delta(t)$ of the equilibrium point $x=0,$ and we came to the surprising conclusion that the origin may remain \enquote{attractive} even for unbounded (and possibly unknown) perturbations $\delta(t),$ as is demonstrated in Example~\ref{example1}. This example shows at the same time that the conditions imposed on the system in Theorem~\ref{theorem_main} cannot be weakened too much. 

Thus the results achieved in this paper contradict the opinion formulated in the classic monograph on dynamical systems \cite[Chapter~9, p.~346]{Khalil}, where it is written: 

\enquote{{\it The origin $x = 0$ may not be an equilibrium point of the perturbed system. We can no longer study stability of the origin as an equilibrium point, nor should we expect the solution of the perturbed system to approach the origin as $t\to\infty.$ The best we can hope for is that $x(t)$ will be ultimately bounded by a small bound, if the perturbation term is small in some sense.}}
%%%%%%%%%%%%%%%%%%%%%%%%%%%%%%%%%%%%%%%%%%%%%%%%%%%
\subsection{Notations} 
Let $\mathbb{R}^n$ denote an $n-$dimensional vector space endowed by any vector norm $\vert\cdot\vert,$ and $\norm{\cdot}$ be an induced norm for matrices, $\norm{A}=\max\{\vert{Ax}\vert;$ $\vert{x}\vert=1\}.$ In the specific situations, when the vector norm is derived from the weighted inner product $(x,y)_P\triangleq y^TPx$ on $\mathbb{R}^n$ and $\vert x\vert_P\triangleq(x,x)_P^{1/2},$ where $P$ is a symmetric and positive definite matrix, we use the notation with the subscript $P,$ $\vert\cdot\vert_P,$ $\norm{\cdot}_P,$ {\it etc.} Obviously, for $P=I$ (the unit matrix on $\mathbb{R}^n$) we obtain the Euclidean norm, $\vert\cdot\vert_I.$ Throughout the whole paper, the superscript \enquote{\,T\,} indicates the transpose operator.
 
We always assume that the function $f:$ $\mathbb{R}^n\times[t_0,\infty)\to\mathbb{R}^n$ is continuously differentiable in $x$ and continuous in $t$ and that perturbation $\delta:$ $[t_0,\infty)\to\mathbb{R}^n$ is continuous.  The perturbing term $\delta(t)$ aggregates all external disturbances which affect the nominal system $\dot x=f(x,t),$ where, as usual, the overdot represents the derivative of the state variable $x=x(t)$ with respect to time $t.$ Let us denote by $J_xf(y,t)$ the Jacobian matrix of $f$ with respect to variable $x$ and evaluated at $(y,t).$   We also assume that the solutions of (\ref{original_system}) are uniquely determined by $x(t_0)$ for all $t\geq t_0.$

For later reference, we introduce two useful relations from the calculus of vector functions.
\begin{lem}\label{integral_eq}
Let the function $f(x,t)$ from $\mathbb{R}^n\times[t_0,\infty)$ to $\mathbb{R}^n$ be a continuously differentiable in $x$ and continuous in $t.$ Then  
\begin{itemize} 
\item[(I)] 
\[
\bigg[\int\limits_0^1 J_xf(\xi x,t)d\xi\bigg]x=f(x,t)-f(0,t);
\]
\end{itemize}
or more generally,
\begin{itemize}
\item[(II)] 
\[
f(x,t)-f(x^*,t)=\bigg[\int\limits_0^1 J_xf(x^*+\xi (x-x^*),t)d\xi\bigg](x-x^*)
\]
\end{itemize}
for all $x,x^*\in\mathbb{R}^n.$
\end{lem}
The proof of Lemma~\ref{integral_eq} is postponed in Appendix and for which we do not claim any originality. 
%%%%%%%%%%%%%
\medskip

\noindent The key role in our analysis plays the logarithmic norm $\mu[A]$ of a matrix $A,$ which is in some sense analogous to a norm, albeit it is not actually a norm in the usual sense, but which gives principally the sharper estimates on asymptotic behavior of the solutions than norms, because $\mu[A]$ may take on also negative values. We define for any real $n\times n$ matrix $A$ the logarithmic norm by the relation
\begin{equation}\label{lognorm_def}
\mu[A]\triangleq\lim\limits_{\theta \to 0^+}\frac{\norm{I+\theta A}-1}{\theta}.
\end{equation}
Specifically, for the Euclidean norm, by \cite{Afanasiev}, \cite{Coppel1}, \cite{Dekker_Verwer}, \cite{Desoer2}, 
\begin{equation}\label{lognorm_euclidean}
\mu_I[A]=\frac12\lambda_{\max}\left({A+A^T}\right),
\end{equation}
where $\lambda_{\max}\left({A+A^T}\right)$ denotes the maximum eigenvalue of the matrix $A + A^T.$ For a general $\mu_P[\cdot]$ see, e.~g. \cite{Hu_Liu},
\[
\mu_P[A]=\frac12\lambda_{\max}\left({\hat A+\hat A^T}\right),\ \hat A=P_0AP_0^{-1}, \ P_0=\sqrt{P}.
\]
The logarithmic norm has the properties \cite{Desoer_Vidyasagar, Desoer2,  Soderlind1, Soderlind2} that are useful in the stability analysis not only for linear systems as we will see later:

\medskip

For any given $n\times n$ real matrices $A, B$ 
\begin{itemize}
\item[(P1)] the limit in (\ref{lognorm_def}) exists; 
\item[(P2)] $\mu[cA+(1-c)B]\leq c\mu[A]+(1-c)\mu[B]$ for all $c\in[0,1]$ (convexity);
\item[(P3)] $\vert \mu[A]-\mu[B]\vert\leq\norm{A-B}$ ($\vert\cdot\vert$ on the left-hand side denotes the absolute value of real number); 
\end{itemize}
\begin{itemize}
\item[(P4)] let $\Phi(t)$ be a fundamental matrix solution for linear time-varying system $\dot x=A(t)x,$ where $A(\cdot):$ $[t_0,\infty)\to\mathbb{R}^{n\times n}$ is a continuous matrix function. Then 
\[
e^{-\int\limits_{\tau}^t \mu[-A(s)]ds}\leq\norm{\Phi(t)\Phi^{-1}(\tau)}\leq e^{\,\int\limits_{\tau}^t \mu[A(s)]ds}
\]
for all $t_0\leq\tau\leq t<\infty;$
\item[(P5)] \cite[p.~34]{Desoer_Vidyasagar} the solution of linear time-varying system $\dot x =A(t)x$ satisfies for all $t\geq t_0$ the inequalities
\[
\vert{x(t_0)}\vert e^{-\int\limits_{t_0}^t \mu[-A(s)]ds}\leq\vert{x(t)}\vert\leq\vert{x(t_0)}\vert e^{\,\int\limits_{t_0}^t \mu[A(s)]ds}.
\]
By the assumption on $A(t)$ and  Property~P3, the integrals above are well-defined because $\mu[A(\cdot)]$ is continuous.
\end{itemize}
\section{Main result}
The main results of the paper are summarized in the following theorem.
\begin{thm}\label{theorem_main}
Let us consider the system (\ref{original_system}),
\[
\dot x=f(x,t)+\delta(t),\quad x\in\mathbb{R}^n, \quad t\geq t_0.
\]
Assume that for some vector norm on $\mathbb{R}^n$ 
\begin{itemize} 
\item[(A1)] there exists a continuous function $\alpha(t)$ and a real constant $\alpha_0>0$ such that 
\[
\mu[J_xf(x,t)]\leq-\alpha(t)\leq-\alpha_0<0\ for\ all\ (x,t)\in\mathbb{R}^n\times[t_0,\infty).
\]
\end{itemize}
Then the difference between any two solutions $x(t)$ and $x^*(t)$ of system (\ref{original_system}) decreases exponentially (and uniformly),
\begin{equation}\label{ineq_main}
\vert x(t)-x^*(t)\vert\leq e^{-\alpha_0(t-t_0)}\vert x(t_0)-x^*(t_0)\vert,\quad t\geq t_0,
\end{equation}
that is, the system (\ref{original_system}) is GIS in the sense of Definition~\ref{GIS}. 

\medskip

In addition, 
\begin{itemize}
\item[(A2)] if the ratio $\frac{\vert f(0,t)+\delta(t)\vert}{\alpha(t)}\to 0$ as $t\to\infty,$ 
\end{itemize}
then all solutions of (\ref{original_system}) converge to $0$ as $t\to\infty.$
\end{thm}
\begin{pf} 
First we prove that the inequality (\ref{ineq_main}) holds. Let us denote by $z$ the difference $z(t)=x(t)-x^*(t),$ $t\geq t_0.$ Observe that $z$ is equal to the solution of linear time-varying system
\[
\dot z=A(t)z,\ z(t_0)=x(t_0)-x^*(t_0)
\]
where, by Lemma~\ref{integral_eq},
\[
A(t)\triangleq \int\limits_0^1 J_xf(x^*(t)+\xi (x(t)-x^*(t)),t)d\xi.
\]
Due to the convexity of the logarithmic norm, Jensen inequality and by the assumption, we obtain
\[
\mu[A(t)]=\mu\bigg[\int\limits_0^1 J_xf(x^*(t)+\xi (x(t)-x^*(t)),t)d\xi\bigg]
\]
\[
\leq\int\limits_0^1 \mu\bigg[J_xf(x^*(t)+\xi (x(t)-x^*(t)),t)\bigg]d\xi\leq-\int\limits_0^1 \alpha(t)d\xi=-\alpha(t)\leq-\alpha_0.
\]
Applying the consistency of the operator norm with the vector norm that induces it and Property~P4 of the logarithmic norm to $z(t)=\Phi(t)\Phi^{-1}(t_0)z(t_0)$ we get (\ref{ineq_main}).

\medskip

Now we prove the second part of Theorem~\ref{theorem_main}, the eventual convergence of all solution to $0$ as $t\to\infty.$
Observe that the solution $x(\cdot)$ of (\ref{original_system}) is equal to the solution of the linear time-varying system
\[
\dot x = \tilde A(t)x+f(0,t)+\delta(t),
\]
where
\[
\tilde A(t)=\int\limits_0^1 J_xf(\xi x(t),t)d\xi.
\]
By similar argument as above, 
\[
\mu[\tilde A(t)]=\mu\bigg[\int\limits_0^1 J_xf(\xi x(t),t)d\xi\bigg]
\]
\begin{equation}\label{est_tildeA}
\leq\int\limits_0^1 \mu\bigg[J_xf(\xi x(t),t)\bigg]d\xi\leq-\int\limits_0^1 \alpha(t)d\xi=-\alpha(t)\leq-\alpha_0.
\end{equation}
Using the variation constant formula, we get
\begin{equation*}
x(t)=\tilde\Phi(t)\bigg[\tilde\Phi^{-1}(t_0)x(t_0)+\int\limits_{t_0}^t \tilde\Phi^{-1}(\tau)[f(0,\tau)+\delta(\tau)]d\tau\bigg],
\end{equation*}
that is,
\begin{equation*}
\vert{x(t)}\vert\leq \vert{x(t_0)}\vert e^{\,\int\limits_{t_0}^t\mu[\tilde A(s)]ds}+\int\limits_{t_0}^t e^{\,\int\limits_{\tau}^t\mu[\tilde A(s)]ds}\vert f(0,\tau)+\delta(\tau)\vert d\tau.
\end{equation*}
Obviously, by Assumption~A1, $\vert{x(t_0)}\vert e^{\,\int\limits_{t_0}^t\mu[\tilde A(s)]ds}\to 0$ (exponentially) for $t\to\infty$ and so it remains to analyze the second term on the right-hand side of the above inequality. We have,
\[
\int\limits_{t_0}^t e^{\,\int\limits_{\tau}^t\mu[\tilde A(s)]ds}\vert f(0,\tau)+\delta(\tau)\vert d\tau= e^{\,\int\limits_{t_0}^t\mu[\tilde A(s)]ds}\int\limits_{t_0}^t e^{-\int\limits_{t_0}^{\tau}\mu[\tilde A(s)]ds}\vert f(0,\tau)+\delta(\tau)\vert d\tau
\]
\begin{equation}\label{estimate_for_limit}
=\frac{\int\limits_{t_0}^t e^{-\int\limits_{t_0}^{\tau}\mu[\tilde A(s)]ds}\vert f(0,\tau)+\delta(\tau)\vert d\tau}{e^{-\int\limits_{t_0}^t\mu[\tilde A(s)]ds}},
\end{equation}
and the L'Hospital rule yields
\[
\lim\limits_{t\to\infty}\frac{\frac{d}{dt}\int\limits_{t_0}^t e^{-\int\limits_{t_0}^{\tau}\mu[\tilde A(s)]ds}\vert f(0,\tau)+\delta(\tau)\vert d\tau}{\frac{d}{dt}e^{-\int\limits_{t_0}^t\mu[\tilde A(s)]ds}}
\]
\[
=\lim\limits_{t\to\infty}\frac{e^{-\int\limits_{t_0}^t\mu[\tilde A(s)]ds}\vert f(0,t)+\delta(t)\vert }{e^{-\int\limits_{t_0}^t\mu[\tilde A(s)]ds}(-\mu[\tilde A(t)])}=\lim\limits_{t\to\infty}\frac{\vert f(0,t)+\delta(t)\vert }{-\mu[\tilde A(t)]}.
\]
Now, from (\ref{est_tildeA}), $-\mu[\tilde A(t)]\geq\alpha(t)$ which, together with Assumption~A2, gives the statement of the second part of Theorem~\ref{theorem_main}.
\end{pf}
\begin{rmk}\label{Demidovich1}
A great Russian mathematician and one of the pioneers in the area of stability of dynamical systems, B.P.~Demidovich showed, see e.~g. \cite{Pavlov} or the original source in Russian \cite{Demidovich}, that if, for some positive definite matrix $P =P^T>0$, the matrix
\begin{equation}\label{Demidovich}
J(x,t)=\frac12\left[PJ_xF(x,t)+J^T_xF(x,t)P\right]
\end{equation} 
is negative definite uniformly in $(x,t)\in\mathbb{R}^n\times\mathbb{R}$ then for any two solutions $x(t)$ and $x^*(t)$ of the dynamical system $\dot x=F(x,t)$ is
\[ 
\vert x(t)-x^*(t)\vert_I\leq K e^{-\alpha(t-t_0)}\vert x(t_0)-x^*(t_0)\vert_I
\]
for all $t\geq t_0$ and some independent on $x$ and $x^*$ constants $K,\alpha>0.$ However, this condition is not very-well suited for reasoning about the convergence of all solutions to $0$ as $t\to\infty$ if \[
(F(0,t)=)\,f(0,t)+\delta(t)\neq0
\]
because we cannot set $x^*(t)\equiv0.$ 

In the context of logarithmic norm, Demidovich condition (\ref{Demidovich}) is equivalent to the existence of  positive definite symmetric matrix $P$ such that $\mu_P[J_xF(x,t)]\leq-\alpha<0.$  

In fact, as follows from \cite{Dekker_Verwer} and \cite{Hu_Liu},
\[
\mu_P[A]=\max\limits_{x\neq 0}\frac{(Ax,x)_P}{\vert x\vert^2_P}=\max\limits_{x\neq 0}\frac{x^T(PA+A^TP)x}{2\vert x\vert^2_P},\ A=J_xF(x,t).
\]
But, taking into account that not every norm comes from an inner product, our result strengthens the Demidovich results.
\end{rmk}
\begin{rmk}\label{Demidovich2}
The condition in Assumption~A1 might be relaxed to 
\[
\forall (x,t)\in\mathbb{R}^n\times[t_0,\infty):\  \mu[J_xf(x,t)]\leq-\alpha(t), \quad \int\limits_{t_0}^{\infty}\alpha(\tau)d\tau=\infty
\]
to obtain only asymptotic stability of solutions (not uniform and not exponential, in general),
\[
\vert x(t)-x^*(t)\vert\leq e^{-\int\limits_{t_0}^{t}\alpha(\tau)d\tau} \vert x(t_0)-x^*(t_0)\vert, \quad t\geq t_0.
\]
Recall that proof by L'Hospital rule requires $\alpha(t)>0$ in some left neighborhood of $t=\infty.$ 

Notice also that, albeit under these circumstances the system may not satisfy the conditions for GIS from Definition~\ref{GIS},  still all solutions converge to one another as $t\to\infty.$ Thus we have extended the results presented in \cite{Lohmiller} to more general type of convergence and also to potentially unbounded perturbation $\delta(t)$ of the nominal system $\dot x=f(x,t).$
\end{rmk}
\section{Simulation experiments}
\begin{figure}[H] 
\captionsetup{singlelinecheck=off} 
   \centerline{
    \hbox{
     \psfig{file=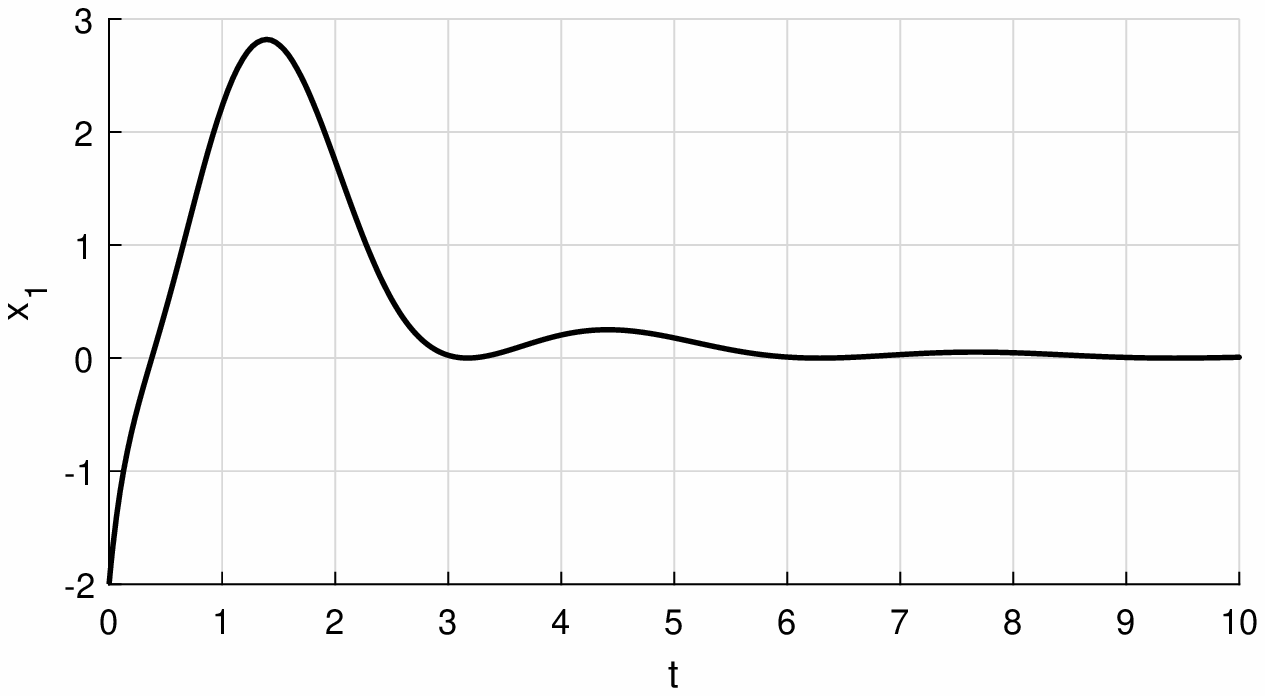,width=5.0cm, clip=}
     \hspace{1.cm}
     \psfig{file=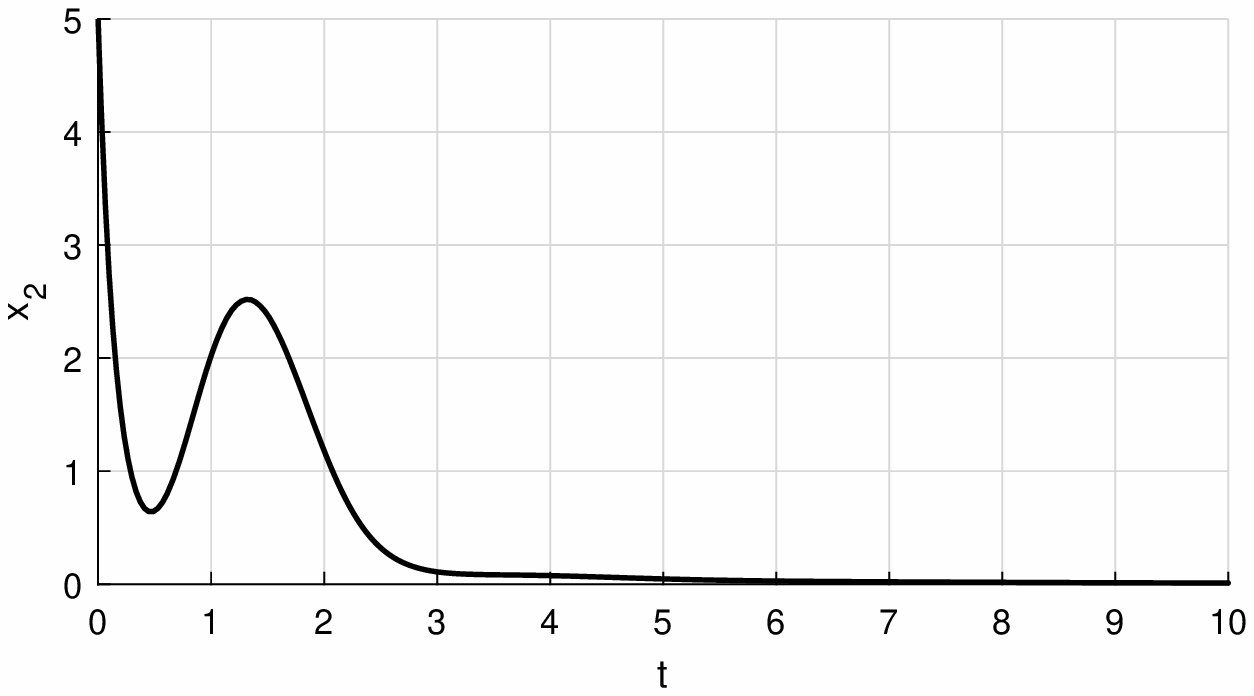,width=5.0cm,clip=}
    }
   }  
\caption{The numerical solution $x(t)=\left(x_1(t), x_2(t)\right)^T$ of the system $\dot x =f(x,t)+\delta(t),$ where $f(x,t)$ is given by (\ref{eq:example1}) with $b=5,$ $\phi(t)=-6-t^3,$ the admissible (unbounded) perturbation  $\delta(t)=\big(5\sin^2\left(t\right),\, t\big)^T$  and the initial state $x(0)=(-2, \ 5)^T.$ }
\label{solution_example1}
\end{figure}

\begin{figure}[H]  
\captionsetup{singlelinecheck=off} 
   \centerline{
    \hbox{
     \psfig{file=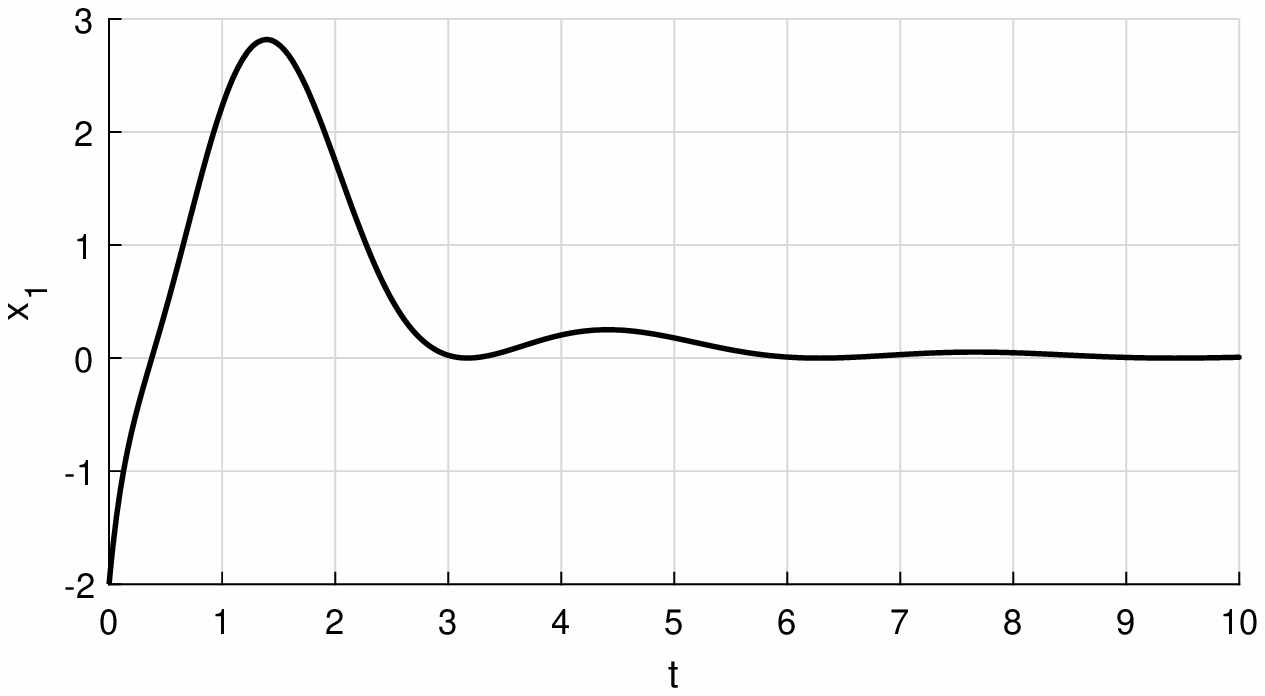,width=5.0cm, clip=}
     \hspace{1.cm}
     \psfig{file=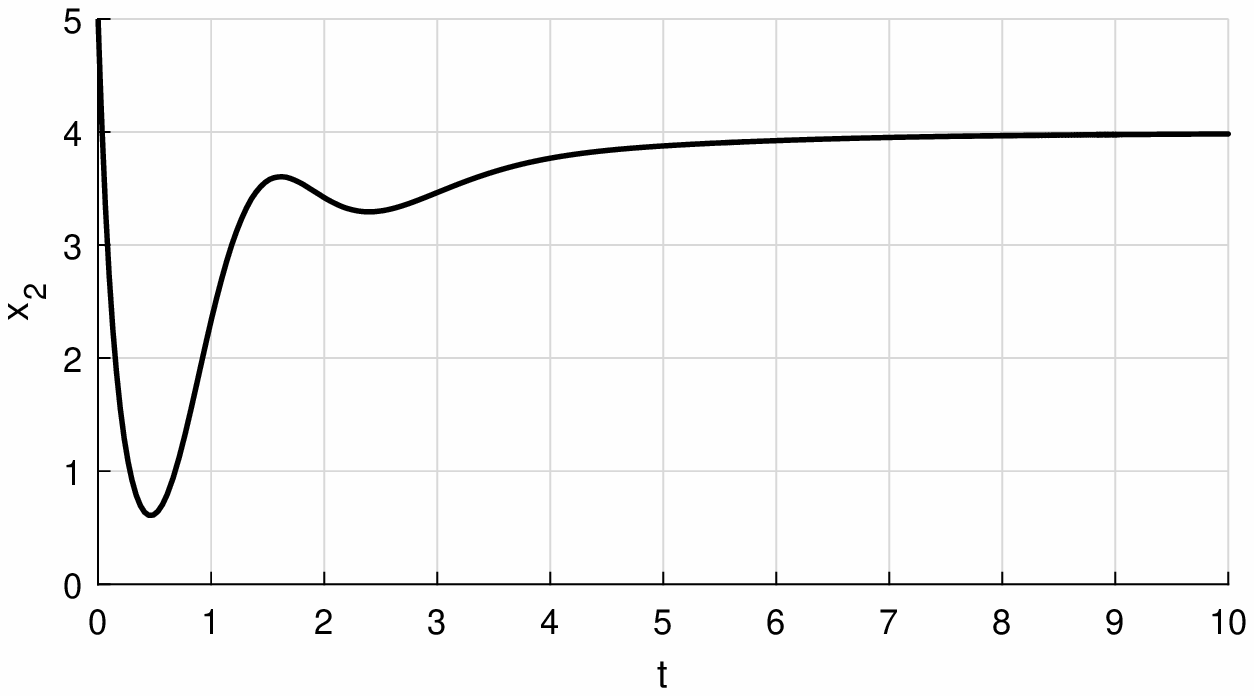,width=5.0cm,clip=}
    }
   }  
\caption{The numerical solution $x(t)=\left(x_1(t), x_2(t)\right)^T$ of the system $\dot x =f(x,t)+\delta(t),$ where $f(x,t)$ is given by (\ref{eq:example1}) with $b=5,$ $\phi(t)=-6-t^3,$  the borderline perturbation $\delta(t)=\big(5\sin^2\left(t\right),\, 4t^3\big)^T$  and the initial state $x(0)=(-2, \ 5)^T.$ }
\label{solution_example1b}
\end{figure}

\begin{ex}\label{example1}
As an academic example, let us consider the planar nonlinear system $\dot x =f(x,t)+\delta(t),$ $t\geq t_0$ with
\begin{equation}\label{eq:example1}
f(x,t)=\big(\phi(t)x_1 +\sin\left(x_1\right),\ bx_1+ [2+\phi(t)]x_2+\sin\left(x_2\right)\big)^T,
\end{equation}
where $\phi(t)$ is an arbitrary scalar continuous function on $[t_0,\infty)$ and $b$ is a real constant. By (\ref{lognorm_euclidean}) and with the help of MATLAB code, we have 
\[
\mu_I\big[J_x f(x,t)\big]=\phi(t)+\frac12\big[\cos\left(x_{1}\right)+\cos\left(x_{2}\right)+\sqrt{\vartheta }\big]+1,
\] 
where
\[
\vartheta = b^2+[{\cos\left(x_{1}\right)}-{\cos\left(x_{2}\right)}]^2-4\,\cos\left(x_{1}\right)+4\,\cos\left(x_{2}\right)+4
\]
\[
=b^2+ [{\cos\left(x_{1}\right)}-{\cos\left(x_{2}\right)}-2]^2\geq 0.
\]
For example, if we choose $b=5$ and $\phi(t)=-6-t^3,$ the Assumptions~A1 and A2 of Theorem~\ref{theorem_main} hold for $\alpha(t)=0.5+t^3,$ $t\geq 0(=t_0)$ and so the vanishing of all solutions as $t\to\infty$ is ensured for perturbations satisfying $\vert\delta(t)\vert_I=o(t^3)$ in Landau's little-o notation.  It means, that the system is GIS in the sense of Definition~\ref{GIS} and, in addition, all solutions converge to $0$ as long as the perturbing term $\delta(t)$ (its $\vert\cdot\vert_I-$norm, to be more precise) is of the order less than $t^3$ as $t\to\infty,$ demonstrating the global robust stability of the equilibrium point $x=0$ of the nominal system ($\delta=0$) even for unbounded perturbations $\delta.$ The results of simulation experiments are shown in Fig.~\ref{solution_example1} and Fig.~\ref{solution_example1b}, where for the simulation purpose we selected one representative from the class of admissible perturbations (Fig.~\ref{solution_example1}) and the borderline case for the second simulation experiment. The dynamics of the system on Fig.~\ref{solution_example1b} indicates that Assumption~A2, ensuring the convergence to zero  of all solutions as $t\to\infty,$ cannot be weakened too much. 

The limiting value $(0,\, 4)^T$ can be in this particular example calculated explicitly thinking as follows: Separately analyzing the first equation by using Theorem~\ref{theorem_main} for $\delta=(\delta_1,\, \delta_2)^T=(5\sin^2\left(t\right),\, 4t^3)^T$ and $\phi(t)=-6-t^3,$ we obtain that $x_1(t)\to0$ as $t\to\infty.$ In the second equation, transforming the second component of state vector by $x_2=\tilde x_2+4$ and identifying  $(bx_1(t)-16)$ as an inhomogeneous term $\tilde\delta_2(t),$ we get
the scalar differential equation 
\[
\dot{\tilde x}_2=\tilde f_2(\tilde x_2,t)+\tilde\delta_2(t),
\]
where
\[ 
\tilde f_2(\tilde x_2,t)\triangleq f_2(x_1,\tilde x_2+4,t)-bx_1+\delta_2(t)=-(4+t^3)\tilde x_2+\sin(\tilde x_2+4).
\]
Then Theorem~\ref{theorem_main} implies that $\tilde x_2(t)\to0$ as $t\to\infty$ because 
\[
\frac{\vert\tilde f_2(0,t)+\tilde\delta_2(t)\vert}{\tilde\alpha(t)}=\frac{\vert \sin(4)+bx_1(t)-16\vert}{3+t^3} \to 0\ \mathrm{as}\ t\to\infty,
\]
which is what we have to prove.
\end{ex} 
\section*{Conclusions}
In this paper, the new result for assessment of the global incremental stability of the nonlinear systems $\dot x=f(x,t)+\delta(t)$ is derived.  Roughly speaking, we have established here the sufficient condition for convergence of any two solutions of a system to each other and another condition for convergence of all solutions to the origin $x=0,$ which may or may not be the equilibrium position for the nominal system $\dot x=f(x,t).$

The fundamental advantage of the used approach based on the logarithmic norm is the fact that to estimate the norm of transition matrix for auxiliary linear time-varying system associated to the original nonlinear one, we do not need to know the fundamental matrix solution and all necessary estimates are based purely on the linear system's matrix entries. 
\section*{Appendix.}
\noindent For completeness, we provide the proof of Lemma~\ref{integral_eq}.
\begin{proofoflemma2}
We prove Part~I only, in the proof of second statement we proceed analogously. Let $f_i,$ $i=1,\dots,n$ denote the components of $f(x,t)$ and define: $g_i: [0,1]\to\mathbb{R}$ by $g_i(\xi)=f_i(\xi x,t).$ Then we have
\[
f_i(x,t)-f_i(0,t)=g_i(1)-g_i(0)=\int\limits_0^1 g'_i(\xi)d\xi
\]
\[
=\int\limits_0^1\bigg(\sum\limits_{j=1}^n\frac{\partial f_i}{\partial x_j}(\xi x,t)x_j\bigg)d\xi=\sum\limits_{j=1}^n\bigg(\int\limits_0^1 \frac{\partial f_i}{\partial x_j}(\xi x,t) d\xi\bigg)x_j, \ i=1,\dots,n.
\]
Now the statement of lemma follows immediately.
\end{proofoflemma2}

%\section*{References}

%%%%%%%%%%%%%%%%%%%%%%%%%%%% 
\end{document}